\documentclass{amsproc}

\usepackage{amssymb}
\usepackage{graphicx}
\usepackage{amscd}
\usepackage{amsmath}

\theoremstyle{plain}
\newtheorem{theorem}{Theorem}
\newtheorem{lemma}{lemma}

\newtheorem{definition}{Definition}

\newtheorem{remark}{Remark}
\newtheorem{proposition}{Proposition}

\begin{document}

\title[Hilbert spaces of analytic
functions]{Inverse problem for semi-infinite Jacobi matrices and
associated Hilbert spaces of analytic functions}

\author{Alexander Mikhaylov } 
\address{St. Petersburg   Department   of   V.A. Steklov    Institute   of   Mathematics
of   the   Russian   Academy   of   Sciences, 7, Fontanka, 191023
St. Petersburg, Russia and Saint Petersburg State University,
St.Petersburg State University, 7/9 Universitetskaya nab., St.
Petersburg, 199034 Russia.} \email{mikhaylov@pdmi.ras.ru}

\author{Victor Mikhaylov} 
\thanks {  }
\address{St.Petersburg   Department   of   V.A.Steklov    Institute   of   Mathematics
of   the   Russian   Academy   of   Sciences, 7, Fontanka, 191023
St. Petersburg, Russia and Saint Petersburg State University,
St.Petersburg State University, 7/9 Universitetskaya nab., St.
Petersburg, 199034 Russia.} \email{ftvsm78@gmail.com}

\date{July, 2019}

\maketitle

\begin{abstract}
We consider the dynamic problems for the discrete systems with
discrete time associated with finite and semi-infinite Jacobi
matrices. The result of the paper is a procedure of association of
special Hilbert spaces of functions, namely de Branges space,
playing an important role in the inverse spectral theory, with
these systems.
We point out the relationships with the classical moment problems
theory and compare properties of classical Hankel matrices
associated with moment problems with properties of matrices of
connecting operators associated with dynamical systems.
\end{abstract}

\section{Introduction}

For a given sequence of positive numbers $\{a_0, a_1,\ldots\}$ (in
what follows we assume $a_0=1$) and real numbers $\{b_1,
b_2,\ldots \}$, we denote by $A$ the semi-infinite Jacobi matrix
\begin{equation}
\label{Jac_matr}
A=\begin{pmatrix} b_1 & a_1 & 0 & 0 & 0 &\ldots \\
a_1 & b_2 & a_2 & 0 & 0 &\ldots \\
0 & a_2 & b_3 & a_3 & 0 & \ldots \\
\ldots &\ldots  &\ldots &\ldots & \ldots &\ldots
\end{pmatrix}.
\end{equation}
For $N\in \mathbb{N}$, by $A^N$ we denote the $N\times N$ Jacobi
matrix which is a block of (\ref{Jac_matr}) consisting of the
intersection of first $N$ columns with first $N$ rows of $A$.
%
We consider the dynamical system corresponding to a semi-infinite
Jacobi matrix:
\begin{equation}
\label{Jacobi_dyn}
\begin{cases}
u_{n,t+1}\!+\!u_{n,t-1}\!-\!a_nu_{n+1,t}\!-\!a_{n-1}u_{n-1,t}\!-\!b_nu_{n,t}\!\!=0,\,\, t\geqslant 0,\,\,
n\geqslant 1,\\
u_{n,\,-1}=u_{n,\,0}=0,\quad n\geqslant 1 \\
u_{0,\,t}=g_t,\quad t\geqslant 0,
\end{cases}
\end{equation}
where $g=(g_0,g_1,\ldots)$ is a \emph{boundary control}, $g_i\in
\mathbb{C}$, $i=0,1,2,\ldots$, the solution to (\ref{Jacobi_dyn})
is denoted by $u^g.$

De Branges spaces \cite{DBr,DMcK,RR} play an important role in the
inverse spectral theory. In \cite{MM6,MM7} the authors show how to
associate finite-dimensional de Branges spaces with the dynamical
systems of the form (\ref{Jacobi_dyn}). Note that our approach
differs from the classical one and potentially admits the
generalization to the multidimensional systems \cite{MM11}. The
algorithm proposed in \cite{MM7,MM6} is as follows: fixing some
finite time $t=T$ one introduces the \emph{reachable set} of the
dynamical system at this time:
\begin{equation*}
U^T:=\{u^g_{\cdot,T}\,|\, g\in \mathcal{F}^T\}.
\end{equation*}
Then one needs to apply the Fourier transform associated with the
operator corresponding  to the matrix $A$ to elements from $U^T$
and get a linear manifold $\mathcal{F}U^T$. Then this linear
manifold is equipped with the norm defined by the \emph{connecting
operator} $C^T$ associated with the system (\ref{Jacobi_dyn}),
which resulted in the finite-dimensional de Branges space
associated with $A^T$.

Thus for the system (\ref{Jacobi_dyn}), due to the finiteness of
the speed of wave propagation, the described procedure leads to
the finite dimensional Hilbert space of analytic functions
associated with $A^T$ (not with the whole matrix $A$!). The
natural question then is to try to associate some
infinite-dimensional functional spaces with the dynamical system
(\ref{Jacobi_dyn}) with semi-infinite matrix, taking $T\to\infty$.

For a given sequence of \emph{moments} $s_0,s_1,s_2,\ldots$  a
solution of a \emph{Hamburger moment problem} \cite{AH,S} is a
Borel measure $d\rho(\lambda)$ on $\mathbb{R}$ such that
\begin{equation}
\label{MikhaylovAS_Moment_eq} s_k=\int\limits_{-\infty}^\infty
\lambda^k\,d\rho(\lambda),\quad k=0,1,2,\ldots.
\end{equation}
If in a Hamburger moment problem an additional constraint
$\operatorname{supp}d\rho\subset [0,+\infty)$ or
$\operatorname{supp}d\rho\subset [0,1]$  is imposed on the
measure, then such a problem is called \emph{Stieltjes moment
problem} or \emph{Hausdorff moment problem}. The set of solutions
of Hamburger moment problem is denoted by $\mathcal{M}_H$. The
moment problem (moment sequence $s_0, s_1,\ldots$) is called
\emph{determinate} if a solution exists and is unique, if a
solution is not unique, it is called \emph{indeterminate}, the
same notations are used in respect to corresponding measures. The
relationships of classical moments problems and dynamic inverse
problem for the system (\ref{Jacobi_dyn}) are described in
\cite{MM10,MM9}.

The following semi-infinite Hankel matrix is associated with the  moment problem (see \cite{AH})
\begin{equation}
\label{S_matr}
S=\begin{pmatrix} s_0 & s_1 & s_2 & s_3 & \ldots\\
s_1 & s_2 & s_3 & \ldots & \ldots \\
s_2 & s_3 & \ldots & \ldots & \ldots \\
s_3 & \ldots & \ldots & \ldots & \ldots \\
\ldots & \ldots & \ldots & \ldots & \ldots
\end{pmatrix}
\end{equation}
The $N\times N$ block of this matrix is denoted by $S_N$.

The \emph{connecting operator} $C^T$ associated with the system
(\ref{Jacobi_dyn}) (for fixed time $t=T$) plays a central role in
the Boundary Control method \cite{B07,B17}, an approach to inverse
dynamic problems. In \cite{MM9,MM10} the authors shown the simple
relationship between $C^T$ and $S_T$. As $C^T$ was used in the
procedure of the construction of finite-dimensional spaces, the
natural question then is to study properties of ``semi-infinite''
connecting operator $C$ in a way the semi-infinite Hankel matrix
$S$ was studied in \cite{Berg1,Berg3,Y,Y2}.

In the second section we provide all the necessary information on
de Branges spaces. In the third section we list the results for
the finite and semi-infinite operators $S_N,$ $S$ according to
\cite{Berg1,Berg3,Y,Y2}. In the forth section we briefly outline
the results for the dynamical system (\ref{Jacobi_dyn}) and for
dynamical system associated with $A^N$ according to \cite{MM8}. In
the fifth section we describe the procedure of the de Branges
spaces construction in the finite dimensional case according to
\cite{MM6,MM7}. After that in the last section we introduce the
semi-infinite matrix $C$ and compare its properties with ones of
the ``moment problem'' counterpart Hankel matrix $S$. All these
give us the possibility to introduce the infinite-dimensional de
Branges spaces in the limit circle (indeterminate) case for the
system (\ref{Jacobi_dyn}) .

\section{de Branges spaces.}

Here we provide the information on de Branges spaces in accordance
with \cite{R2,RR}. The entire function $E:\mathbb{C}\mapsto
\mathbb{C}$ is called a \emph{Hermite-Biehler function} if
$|E(z)|>|E(\overline z)|$ for $z\in \mathbb{C}_+$. We use the
notation $F^\#(z)=\overline{F(\overline{z})}$. The \emph{Hardy
space} $H_2$ is defined by: $f\in H_2$ if $f$ is holomorphic in
$\mathbb{C}^+$ and
$\sup\limits_{y>0}\int\limits_{-\infty}^\infty|f(x+iy)|^2\,dx<\infty$. Then the
\emph{de Branges space} $B(E)$ consists of entire functions such
that:
\begin{equation*}
B(E):=\left\{F:\mathbb{C}\mapsto \mathbb{C},\,F \text{ entire},
\,\frac{F}{E},\frac{F^\#}{E}\in
H_2\right\}.
\end{equation*}
The space $B(E)$ with the scalar product
\begin{equation*}
[F,G]_{B(E)}=\frac{1}{\pi}\int\limits_{\mathbb{R}}\overline{ F(\lambda)}
{G(\lambda)}\frac{d\lambda}{|E(\lambda)|^2}
\end{equation*}
is a Hilbert space. For any $z\in \mathbb{C}$ the
\emph{reproducing kernel} is introduced by the relation \cite[ Th.
19, p. 50]{DBr}.
\begin{equation}
\label{repr_ker} J_z(\xi):=\frac{\overline{E(z)}E(\xi)-E(\overline
z)\overline{E(\overline \xi)}}{2i(\overline z-\xi)}.
\end{equation}
Then
\begin{equation*}
F(z)=[J_z,F]_{B(E)}=\frac{1}{\pi}\int\limits_{\mathbb{R}}\overline{J_z(\lambda)}
{F(\lambda)}\frac{d\lambda}{|E(\lambda)|^2}.
\end{equation*}
We observe that a Hermite-Biehler function $E(\lambda)$ defines
$J_z$ by (\ref{repr_ker}). The converse is also true \cite[Th. 22,
p 55]{DBr}:
\begin{theorem}
\label{TeorDB} Let $X$ be a Hilbert space of entire functions with
reproducing kernel such that
\begin{itemize}

\item[1)] if $f\in X$ then $f^\#\in X$ and $\|f\|_X=\|f^\#\|_X$,

\item[2)] if $f\in X$ and $\omega\in \mathbb{C}$ such that
$f(\omega)=0$, then $\frac{z-\overline{\omega}}{z-\omega}f(z)\in
X$ and
$\left\|\frac{z-\overline{\omega}}{z-\omega}f(z)\right\|_{X}=\|f\|_{X}$,
\end{itemize}
then $X$ is a de Branges space based on the function
\begin{equation*}
E(z)=\sqrt{\pi}(1-iz)J_i(z)\|J_i\|_X^{-1},
\end{equation*}
where $J_z$ is a reproducing kernel.
\end{theorem}


\section{Classical moment problems, Hankel matrices, minimal eigenvalues and closability. Associated Hilbert space of analytic functions.}

With the semi-infinite matrix $A$ we associate the symmetric
operator $A$ (we keep the same notation) in the space $l_2$,
defined on finite sequences:
\begin{equation}
\label{domA}
D=\left\{\varkappa=(\varkappa_1,\varkappa_2,\ldots)\,|\,\exists\ N\in \mathbb{N}:\   \varkappa_n=0,\,
n\geqslant N  \right\},
\end{equation}
and given by the rule
\begin{align*}
(A\theta)_1&=b_1\theta_1+a_1\theta_2,\\
(A\theta)_n&=a_{n}\theta_{n+1}+a_{n-1}\theta_{n-1}+b_n\theta_n,\quad
n\geqslant 2. 
\end{align*}
By $[\cdot,\cdot]$ we denote the scalar product in $l_2$. For a
given sequence
$\varkappa=\left(\varkappa_1,\varkappa_2,\ldots\right)$ we define
a new sequence
\begin{align*}
&\left(G\varkappa\right)_1=
b_1\varkappa_1+a_1\varkappa_2,\\
&\left(G\varkappa\right)_n=a_{n}\varkappa_{n+1}+a_{n-1}\varkappa_{n-1}+b_n\varkappa_n,\quad
n\geqslant 2.
\end{align*}
The adjoint operator $A^*\varkappa=G\varkappa$ is defined on the
domain
\begin{equation*}
D\left(A^*\right)=\left\{\varkappa=(\varkappa_1,\varkappa_2,\ldots)\in
l_2\,|\, (G\varkappa)\in l_2\right\}.
\end{equation*}
In the limit point case (i.e when $A$ has deficiency indices
$(0,0)$), $A$ is essentially self-adjoint. In the limit circle
case, i.e. when $A$ has deficiency indices $(1,1)$, we denote by
$p(\lambda)=(p_1(\lambda),p_2(\lambda),\ldots)$ and by
$q(\lambda)=(q_1(\lambda),q_2(\lambda),\ldots)$ two solutions of
the difference equation (we set here $a_0$=1):
\begin{equation}
\label{Phi_def}
a_n\phi_{n+1}+a_{n-1}\phi_{n-1}+b_n\phi_n=\lambda\phi_n,\quad
n\geqslant 1,
\end{equation}
satisfying Cauchy data $p_1(\lambda)=1,$
$p_2(\lambda)=\frac{\lambda-b_1}{a_1}$, $q_1(\lambda)=0,$
$q_2(\lambda)=\frac{1}{a_1}$. Thus $p_n(\lambda),\, q_n(\lambda)$
are polynomials of orders $n-1$ and $n-2$. Then \cite[Lemma
6.22]{KSch}
\begin{equation*}
D\left(A^*\right)=D(\overline A)\dot +\mathbb{C}p(0)\dot
+\mathbb{C}q(0),
\end{equation*}
where $\dot +$ denotes the direct sum and $\overline A$ is a
closure of $A$. All self-adjoint extensions of $A$ are
parameterized by $h\in \mathbb{R}\cup\{\infty\}$, are denoted by
$A_{\infty,\,h}$ and are defined on the domain
\begin{equation*}
D(A_{\infty,\,h})=\begin{cases} D(\overline A)\dot +\mathbb{C}(q(0)+hp(0)),\quad h\in \mathbb{R}\\
D(\overline A)\dot +\mathbb{C}p(0),\quad h=\infty.
\end{cases}
\end{equation*}
All the details the reader can find in \cite{S,KSch}. We introduce
the measure
$d\rho_{\infty,\,h}(\lambda)=\left[dE^{A_{\infty,\,h}}_\lambda
e_1,e_1\right]$, where $e_1=(1,0,\ldots)$ and $dE^{A_{\infty,h}}_\lambda$ is the
projection-valued spectral measure of $A_{\infty,\,h}$ such that
$E^{A_{\infty,\,h}}_{\lambda-0}=E^{A_{\infty,\,h}}_{\lambda}$.

Let $\{s_0,s_1.s_2,\ldots\}$ be a moment sequence,  $d\rho\in
\mathcal{M}_H$ be a solution to moment problem
(\ref{MikhaylovAS_Moment_eq}), $S,$ $S_N$ be Hankel matrices
(\ref{S_matr}). The well-known fact on the solvability of the
moments problem \cite{AH,S,KSch} is given in terms of positivity
of these matrices:
\begin{theorem}
\label{ThMomentEx} The moment problem
(\ref{MikhaylovAS_Moment_eq}) has a solution if and only if the
following condition holds:
\begin{equation*}
S_N>0,\quad N=1,2,3,\ldots
\end{equation*}
\end{theorem}
Note that in this case the matrix $A$ is determined by the moments
in a one-to-one manner.

In \cite{Berg1,Berg2,Berg3,Y,Y2} the authors studied properties of
matrix $S_N$ and corresponding semi-infinite matrix $S$. Below we
list the number of important results obtained by them. Having
fixed $N$, we denote by $\lambda_N$ the smallest eigenvalue of
$S_N$:
\begin{equation*}
\lambda_N=\min\{l_k\,|\,l_k\text{ is eigenvalue of $S_N$},\,
k=1,\ldots,N \}.
\end{equation*}

The following theorem was proved in \cite{Berg1}:
\begin{theorem}
\label{BergTH} The moment problem associated with the sequence $\{s_k\}$ is
determined (the matrix $A$ is in the limit point case) if and only
if
\begin{equation*}
\lim_{N\to\infty}\lambda_N\to 0.
\end{equation*}
In the limit circle case
\begin{equation*}
\lim_{N\to\infty}\lambda_N\geqslant
\left(\int_0^{2\pi}l\left(e^{i\theta}\right)\frac{d\theta}{2\pi}\right)^{-1},\quad
l(z)=\left(\sum_{k=0}^\infty|p_k(z)|^2\right)^{-1}.
\end{equation*}
\end{theorem}

The Hankel matrix $S$ give rise to the formally defined
operator $Q$ in $l_2$:
\begin{equation*}
\left(Qf\right)_n=\sum_{m=0}^\infty s_{n+m}f_m,\quad f\in l_2.
\end{equation*}
Then without any a priory assumptions on the sequence $\{s_k\}$,
only the quadratic form of this operator
\begin{equation}\label{quadr_form_S}
S[f,f]=\sum_{m,n\geqslant 0}s_{m+n}\overline {f_m}{f_n}
\end{equation}
is well-defined on $D$. We always assume the positivity condition
(cf. Theorem \ref{ThMomentEx})
\begin{equation*}
\sum_{m,n\geqslant 0}s_{m+n}\overline
{f_m}{f_n}\geqslant 0,\quad f\in D.
\end{equation*}


The following result is obtained in \cite{Berg3}.
\begin{theorem}
If operator $A$ in the limit circle case then operator $S$ defined via quadratic form $(\ref{quadr_form_S})$ is closable.

There exist operator $A$ in the limit point case such that
operator $S$ defined via quadratic form $(\ref{quadr_form_S})$ is
closable.
\end{theorem}
The solutions of the indeterminate moment problem form an infinite
convex set $\mathcal{M}_H$ of measures $M$ with unbounded support,
for which moment identities (\ref{MikhaylovAS_Moment_eq}) hold. In
this set there exist \emph{extremal} measures $\widetilde M\in
\mathcal{M}_H$ such that the set of polynomials
$\mathbb{C}[\lambda]$ is dense in $L_2(\mathbb{R},\widetilde M)$.
These measures correspond to Neumann extension $A_{\infty,h}$ of
operator $A$ were described in the beginning of Section 3, we
denoted them $d\rho_{\infty,h}$, these measures are discrete and
are of the form
\begin{equation*}
\widetilde M=\sum_{k} c_k\delta_{\lambda_k}(\lambda).
\end{equation*}
Then the result of Stieltjes says that if one mass is removed,
then the new measure
\begin{equation*}
M_1:=\widetilde M-c_1\delta_{\lambda_1}(\lambda)
\end{equation*}
is determinate. The existence of examples of the closability of
$S$ in the determinate situation follows from this theorem.

In the indeterminate case the orthonormal with respect to measure
$M$ polynomials form an orthonormal basis in $L_2(\mathbb{R},M)$
if $M$ is extremal. If $M$ is not extremal then they are basis in
$\overline{\mathbb{C}}[\lambda]$ where the closure is assumed in
$L_2(\mathbb{R},M)$. In both cases
$\overline{\mathbb{C}}[\lambda]$ is isomorphic to the space
$\mathcal{E}$ of entire functions of the form
\begin{equation*}
u(z)=\sum_{k=1}^\infty g_kp_k(z), \quad g\in l^2,\quad \text{$p_k$
are orthonormal w.r.t. measure $M$}.
\end{equation*}
Note (see Section 7 of \cite{KSch}) that the reproducing kernel in
this space has a form
\begin{equation}
\label{Kenrnel0} J^\infty_z(\lambda)=\sum_{n=1}^\infty
\overline{p_n(z)}{p_n(\lambda)}
\end{equation}
and the right hand side  converges  uniformly on compact subsets of $\mathbb{C}$ to
holomorphic function on $\mathbb{C}^2$. For this kernel one has
\begin{equation}
\label{Kernel} \int\limits_{-\infty}^\infty
\overline{J^\infty_z}(\lambda){f(\lambda)}dM(\lambda)=f(z),
\end{equation}
for all polynomials $f(z)$.

In the determinate (limit point) case $\mathbb{C}[\lambda]$ is
dense in $L_2(\mathbb{R},M)$, but since the quantity
(\ref{Kenrnel0}) diverges, and as a consequence, the reproducing
kernel is absent, there are no results on the equivalence of
$\mathbb{C}[\lambda]$ to some space of analytic functions.

\section{Dynamical system with discrete time associated with finite and semi-infinite Jacobi matrices.}

Proofs of statements of this section can be found in
\cite{MM8,MM9}.


We consider the dynamical system with discrete time associated
with the matrix $A^N$:
\begin{equation}
\label{Jacobi_dyn_int}
\begin{cases}
v_{n,t+1}\!+\!v_{n,t-1}\!-\!a_nv_{n+1,t}\!-\!a_{n-1}v_{n-1,t}\!-\!b_nv_{n,t}\!=\!0,\ t\geqslant 0,\ 1\leqslant n \leqslant N,\\
v_{n,-1}=v_{n,0}=0,\quad 1\leqslant n\leqslant N+1, \\
v_{0,t}=f_t,\quad v_{N+1,t}=0,\quad t\geqslant 0,
\end{cases}
\end{equation}
where $f=(f_0,f_1,\ldots)$ is a \emph{boundary control}, $f_i\in
\mathbb{C}$, $i=0,1,2,\ldots$. The solution to
(\ref{Jacobi_dyn_int}) is denoted by $v^f$.

The operator corresponding to a finite Jacobi matrix we also
denote by $A^N: \mathbb{R}^N\mapsto \mathbb{R}^N$, it is given by
\begin{equation*}
\begin{cases}
(A\psi)_1=b_1\psi_1+a_1\psi_2,\quad n=1,\\
(A\psi)_n=a_{n}\psi_{n+1}+a_{n-1}\psi_{n-1}+b_n\psi_n,\quad
2\leqslant n\leqslant N-1,\\
(A\psi)_N=a_{N-1}\psi_{N-1}+b_N\psi_N,\quad n=N,
\end{cases}
\end{equation*}

Denote by $\{\lambda_k\}_{k=1}^N$ roots of the equation
$p_{N+1}(\lambda)=0$, it is known \cite{AH} that they are real and
distinct. We introduce the vectors $\phi^n\in \mathbb{R}^N$ by the
rule $\phi^n_i:=p_i(\lambda_n)$, $n,i=1,\ldots,N,$ and define the
numbers $\rho_k$ by
\begin{equation*}
(\phi^k,\phi^l)=\delta_{kl}\rho_k,\quad k,l=1,\ldots,N,
\end{equation*}
where $(\cdot,\cdot)$ is a scalar product in $\mathbb{R}^N$, and $\delta_{kl}$ is a Kronecker  symbol.
\begin{definition}
The set of pairs
\begin{equation*}
\{\lambda_k,\rho_k\}_{k=1}^N
\end{equation*}
is called spectral data of the operator $A^N$.
\end{definition}

The spectral function of operator $A^N$ is introduced by the rule
\begin{equation*}
\rho_{N}(\lambda):=\sum_{\{k\,|\,\lambda_k<\lambda\}}\frac{1}{\rho_k}.
\end{equation*}

The results of \cite[Sections 4.5, 5.5]{AT}  imply that in the
limit circle case $d\rho_N\to d\rho_{\infty,\,h}$ $*-$weakly as
$N\to\infty,$ where
\begin{equation}
\label{alpha_def} h=-\lim_{n\to\infty}\frac{q_n(0)}{p_n(0)}.
\end{equation}

The outer space of dynamical systems (\ref{Jacobi_dyn_int}),
(\ref{Jacobi_dyn}) is $\mathcal{F}^T:=\mathbb{C}^T$,
$\mathcal{F}^T\ni g,f=(f_0,f_1,\ldots,f_{T-1})$ with the inner
product
$(f,g)_{\mathcal{F}^T}:=\sum_{n=0}^{T-1}\overline{f_i}{g_i}$.

The input $\longmapsto$ output correspondences in systems
(\ref{Jacobi_dyn_int}), (\ref{Jacobi_dyn}) are realized by
\emph{response operators}: $R^T_N,\,R^T:\mathcal{F}^T\mapsto
\mathbb{C}^T$ defined by rules
\begin{eqnarray*}
\left(R^T_Nf\right)_t=v^f_{1,\,t}=\left(r^N*f_{\cdot-1}\right)_t=\sum_{s=0}^{t}r^N_sf_{t-1-s},
\quad t=1,\ldots,T,\\
\left(R^Tf\right)_t=u^f_{1,\,t}=\left(r*f_{\cdot-1}\right)_t=\sum_{s=0}^{t}r_sf_{t-1-s},
\quad t=1,\ldots,T,
\end{eqnarray*}
where $r^N=(r_0^N,r_1^N,\ldots,r_{T-1}^N)$,
$r=(r_0,r_1,\ldots,r_{T-1})$ are \emph{response vectors},
convolution kernels of response operators.

Let $\mathcal{T}_k(2\lambda)$ be Chebyshev polynomials of the
second kind, i.e. they satisfy
\begin{equation*}
\begin{cases}
\mathcal{T}_{t+1}(\lambda)+\mathcal{T}_{t-1}(\lambda)-\lambda \mathcal{T}_{t}(\lambda)=0,\\
\mathcal{T}_{0}(\lambda)=0,\,\, \mathcal{T}_1(\lambda)=1.
\end{cases}
\end{equation*}

\begin{proposition}
The solution $v^f$ to system $(\ref{Jacobi_dyn_int})$ and the response vector $r^N$
admit representations
\begin{eqnarray}
&v^f_{n,t}=\int\limits_{-\infty}^\infty \sum\limits_{k=1}^t
\mathcal{T}_k(\lambda)f_{t-k}\phi_n(\lambda)\,d\rho^{N}(\lambda), \quad t\in \mathbb{N},\,n=1,\ldots N,\label{Jac_sol_spectr}\\
&r_{t-1}^{N}=\int\limits_{-\infty}^\infty
\mathcal{T}_t(\lambda)\,d\rho^{N}(\lambda),\quad t\in \mathbb{N}.
\label{Jac_resp_spectr}
\end{eqnarray}
\end{proposition}
\begin{remark}
The solution $u^f$  and response vector $r$ corresponding to
system $(\ref{Jacobi_dyn})$ with semi-infinite Jacobi matrix $A$
admit representations $(\ref{Jac_sol_spectr})$,
$(\ref{Jac_resp_spectr})$ with $d\rho^{N}$ substituted by any
$d\rho(\lambda)\in \mathcal{M}_H$ and $n\in \mathbb{N}$.
\end{remark}

The \emph{inner space} of dynamical system (\ref{Jacobi_dyn_int})
is $\mathcal{H}^N:=\mathbb{C}^N$, $h\in \mathcal{H}^N$,
$h=(h_1,\ldots, h_N)$, $h,m\in \mathcal{H}^N$,
$(h,m)_{\mathcal{H}^N}:=\sum_{k=1}^N\overline{h_k}{m_k}$,
$v^f_{\cdot,\,T}\in \mathcal{H}^N$ for all $T\in \mathbb{N}$. For
the system (\ref{Jacobi_dyn_int}) the \emph{control operator}
$W^T_{N}:\mathcal{F}^T\mapsto \mathcal{H}^N$ is defined by the
rule
\begin{equation*}
W^T_{N}f:=\left(v^f_{1,\,T},\ldots,v^f_{N,\,T}\right).
\end{equation*}
The set
\begin{equation*}
\mathcal{U}^T:=W_N^T \mathcal{F}^T=\left\{\left(v^f_{1,\,T},\ldots,v^f_{N,\,T}\right)\,\bigl|\,
f\in \mathcal{F}^T\right\}
\end{equation*}
is called reachable. For the system (\ref{Jacobi_dyn}) the control operator
$W^T:\mathcal{F}^T\mapsto \mathcal{H}^T$ is introduced by
\begin{equation*}
W^Tf:=\left(u^f_{1,\,T},\ldots,u^f_{T,\,T}\right).
\end{equation*}

This operator admits  the representation $W^T= W_TJ_T$:
\begin{equation}
\label{WT_1} W_T\!\!=\!\!\begin{pmatrix}
a_0 & w_{1,1} & w_{1,2} & \ldots & w_{1,T-1}\\
0 & a_0a_1 & w_{2,2} &  \ldots & w_{2, T-1}\\
\cdot & \cdot & \cdot & \cdot & \cdot \\
0 & \ldots & \prod_{j=1}^{k-1}a_j& \ldots & w_{k, T-1}\\
\cdot & \cdot & \cdot  & \cdot & \cdot \\
0 & 0 & 0  & \ldots & \prod_{j=1}^{T-1}a_j
\end{pmatrix},\ J_T\!\!=\!\!
\begin{pmatrix}
0 & 0 & 0 & \ldots & 1\\
0 & 0 & 0 & \ldots  & 0\\
\cdot & \cdot & \cdot & \cdot &  \cdot \\
0 & \ldots & 1& 0  & 0\\
\cdot & \cdot & \cdot & \cdot &  \cdot \\
1 & 0 & 0 & 0 &  0
\end{pmatrix}.
\end{equation}

Everywhere below we substantially use the finiteness of the speed
of wave propagation in systems (\ref{Jacobi_dyn_int}),
(\ref{Jacobi_dyn}), which implies the following remark.
\begin{remark}
\label{Rem1}
Solution $v^f$ to system $(\ref{Jacobi_dyn})$ and solution $u^f$ to system $(\ref{Jacobi_dyn_int})$ satisfy
\begin{eqnarray*}
&u^f_{n,\,t}=v^f_{n,\,t},\quad n\leqslant t\leqslant N.\notag\\
& 
\quad W^N=W^N_{N}, \quad
r_t=r^N_t,\quad t=0,\ldots,2N-1. 
\end{eqnarray*}
\end{remark}

The \emph{connecting operator} for the system
(\ref{Jacobi_dyn_int}) $C^T_{N}: \mathcal{F}^T\mapsto
\mathcal{F}^T$ is defined via the quadratic form: for arbitrary
$f,g\in \mathcal{F}^T$ we set
\begin{equation*}
\left(C^T_{N} f,g\right)_{\mathcal{F}^T}:=\left(W^T_{N}f,W^T_{N}g\right)_{\mathcal{H}^N}=\left(v^f_{\cdot,\,T},
v^g_{\cdot,\,T}\right)_{\mathcal{H}^N}.
\end{equation*}
For the system (\ref{Jacobi_dyn}) the connecting operator $C^T:
\mathcal{F}^T\mapsto \mathcal{F}^T$ is introduced by the rule:
\begin{equation}
\label{CT_def}
\left(C^T f,g\right)_{\mathcal{F}^T}:=\left(W^Tf,W^Tg\right)_{\mathcal{H}^T}=\left(u^f_{\cdot,\,T},
u^g_{\cdot,\,T}\right)_{\mathcal{H}^T}.
\end{equation}

In \cite{MM7,MM8} the following formulas were obtained:
\begin{proposition}
The matrix of the connecting operator $C_N^T$  for systems $(\ref{Jacobi_dyn_int})$ and the matrix of the connecting operator $C^T$ for system $(\ref{Jacobi_dyn})$ admit spectral representations
\begin{eqnarray*}
\{C^T_{N}\}_{l+1,\,m+1}=\int\limits_{-\infty}^\infty
\mathcal{T}_{T-l}(\lambda)\mathcal{T}_{T-m}(\lambda)\,d\rho^{N}(\lambda),
\quad l,m=0,\ldots,T-1, \label{MikhaylovAS_SP_mes_d}\\
\{C^T\}_{l+1,\,m+1}=\int\limits_{-\infty}^\infty
\mathcal{T}_{T-l}(\lambda)\mathcal{T}_{T-m}(\lambda)\,d\rho(\lambda),
\quad l,m=0,\ldots,T-1,\notag
\end{eqnarray*}
in the last equality we can take any $d\rho(\lambda)\in
\mathcal{M}_H$. The following dynamic
representation valid if $T\leqslant N$: 
    \begin{eqnarray}
&C^T=C^T_{N},\quad \{C^T\}_{ij}=\sum\limits_{k=0}^{T-\max{\{i,j\}}}r_{|i-j|+2k},\nonumber\\
&C^T=
\begin{pmatrix}
r_0\!+\!r_2\!+\!\ldots\!+\!r_{2T-2} & \cdot&\ldots        &  r_T\!+\!r_{T-2}   &  r_{T-1} \\
r_1\!+\!r_3\!+\!\ldots\!+\!r_{2T-3} & \cdot&\ldots        &r_{T-1}\!+\!r_{T-3} &  r_{T-2} \\
\cdot       & \cdot&\cdot         & \cdot          & \cdot    \\
r_{T-3}\!+\!r_{T-1}\!+\!r_{T+1} & \cdot&r_0\!+\!r_2\!+\!r_4   & r_1\!+\!r_3        & r_2      \\
r_{T}\!+\!r_{T-2}           & \cdot&r_1\!+\!r_3       & r_0\!+\!r_2        & r_1       \\
r_{T-1}                 & \cdot&r_2           & r_1            &r_0
\end{pmatrix}\label{MikhaylovAS_SP_mes_d1}
\end{eqnarray}
\end{proposition}

We introduce matrices $C_T:=J_TC^TJ_T$ thus this matrix keeps the
structure of $C^T$ but ``is filled''  from the upper left corner.
Then we have
\begin{equation*}
C_T=\left(W_T\right)^* W_T.
\end{equation*}

\section{De Branges spaces for finite Jacobi matrices}

By $d\rho$ we denote the spectral measure of $A$ in the limit
point case or $d\rho_{\infty,h}$ with $h$ defined in
(\ref{alpha_def}) in the limit circle case. According to \cite{AT}
this measure give rise to the Fourier transform $F: l^2\mapsto
L_2(\mathbb{R},d\rho)$, defined by the rule:
\begin{equation*}
(Fa)(\lambda)=\sum_{n=0}^\infty a_k p_k(\lambda),\quad
a=(a_0,a_1,\ldots)\in l^2,
\end{equation*}
where $p_k(\lambda)$ is a solution to (\ref{Phi_def}) satisfying
Cauchy data $p_1(\lambda)=1,$
$p_2(\lambda)=\frac{\lambda-b_1}{a_1}$. The inverse transform and
Parseval identity have forms:
\begin{eqnarray}
a_k=\int\limits_{-\infty}^\infty
(Fa)(\lambda)p_k(\lambda)\,d\rho(\lambda),\notag\\
\sum_{k=0}^\infty \overline{a_k}{b_k}=\int\limits_{-\infty}^\infty
\overline{(Fa)(\lambda)}{(Fb)(\lambda)}\,d\rho(\lambda).\label{JM_parseval}
\end{eqnarray}
We assume that $T$ is fixed and $f\in \mathcal{F}^T$. Then for
such  control and for $\lambda\in \mathbb{C}$ we have the
following representation \cite{MM7} for the Fourier transform of
the solution to (\ref{Jacobi_dyn}) at $t=T$:
\begin{equation*}
\left(Fu^f_{\cdot,T}\right)(\lambda)=\sum_{k=1}^T
\mathcal{T}_k(\lambda)f_{T-k},\quad \lambda\in \mathbb{C}.
\end{equation*}
Now we 
introduce the linear manifold of Fourier images of states of
dynamical system (\ref{Jacobi_dyn}) at time $t=T$, i.e. the
Fourier image of the reachable set:
\begin{equation}
\label{De_Br_sp}
B_{A}^T\!\!:=\!F\mathcal{U}^T\!\!=\!\left\{\left(Fu^{J_Tf}_{\cdot,T}\right)(\lambda)\,|\,
J_T f\in \mathcal{F}^T\right\}=\left\{\sum_{k=1}^T
\mathcal{T}_k(\lambda)f_{k}\,|\, f\in \mathcal{F}^T\right\}.
\end{equation}
Note that $B^T_A=\mathbb{C}_{T-1}[\lambda]$. It would be
preferable for us to use $C_T$ instead of $C^T$, although that
leads to changing of some formulas comparing to \cite{MM7,MM10}.

We equip $B_{A}^T$ with the scalar product defined by the rule:
\begin{gather}
[F,G]_{B^T_{A}}:=\left(C_Tf,g\right)_{\mathcal{F}^t}, \quad F,G\in
B^T_{A},\label{JM_scalprod}\\
F(\lambda)=\sum_{k=1}^T
\mathcal{T}_k(\lambda)f_{k},\,G(\lambda)=\sum_{k=1}^T
\mathcal{T}_k(\lambda)g_{k},\,f,g\in \mathcal{F}^T.\notag
\end{gather}
Evaluating (\ref{JM_scalprod}) making use of (\ref{JM_parseval})
yields:
\begin{multline}
[F,G]_{B^T_{A}}=\left(C_Tf,g\right)_{\mathcal{F}^T}=\left(C^TJ_Tf,J_Tg\right)_{\mathcal{F}^T}
=\left(u^{J_Tf}_{\cdot,T},u^{J_Tg}_{\cdot,T}\right)_{\mathcal{H}^T}\label{DBr_Scal}\\
=\int\limits_{-\infty}^\infty
\overline{(Fu^{J_Nf}_{\cdot,N})(\lambda)}{(Fu^{J_Ng}_{\cdot,N})(\lambda)}\,d\rho_{\infty,\,h}(\lambda)
=\int\limits_{-\infty}^\infty
\overline{F(\lambda)}{G(\lambda)}\,d\rho_{\infty,\,h}(\lambda).
\end{multline}
Where the last equality is due to the finite speed of wave
propagation in (\ref{Jacobi_dyn}), (\ref{Jacobi_dyn_int}). In
\cite{MM8} the authors proved the following
\begin{theorem}
The vector $(r_0,r_1,r_2,\ldots,r_{2N-2})$ is a response vector
for the dynamical system (\ref{Jacobi_dyn_int}) if and only if the
matrix of connecting operator $C^N$  defined by $(\ref{CT_def})$,
$(\ref{MikhaylovAS_SP_mes_d1})$ with $T=N$ is positive definite.
\end{theorem}
This theorem shows that (\ref{DBr_Scal}) is a scalar product in
$B^T_{A}$. But we can say even more \cite{MM7}:
\begin{theorem}
By dynamical system with discrete time  $(\ref{Jacobi_dyn})$ one
can construct the de Branges space by (\ref{De_Br_sp})
As a set of functions it coincides with the space of Fourier
images of states of dynamical system $(\ref{Jacobi_dyn})$ at time
$t=T$ $($or what is the same, states of $(\ref{Jacobi_dyn_int})$
with $N=T$ at the same time$)$, i.e. the Fourier image of a
reachable set, and is the set of polynomials with real
coefficients of the order less or equal to $T-1$. The norm in
$B_{A}^T$ is defined via the connecting operator:
\begin{equation*}
[F,G]_{B^T_{A}}:=\left(C_Tf,g\right)_{\mathcal{F}^T},\quad F,G\in
B^T_{A},
\end{equation*}
where
\begin{equation*}
F(\lambda)=\sum_{k=1}^T
\mathcal{T}_k(\lambda)f_{k},\,G(\lambda)=\sum_{k=1}^T
\mathcal{T}_k(\lambda)g_{k},\,f,g\in \mathcal{F}^T.
\end{equation*}
The reproducing kernel has a form
\begin{equation*}
J^T_z(\lambda)=\sum_{k=1}^T \mathcal{T}_k(\lambda)(j^z_T)_{k},
\end{equation*}
where $j^z_T$ is a solution to Krein-type equation
\begin{equation}
\label{Krein_eq}
C_Tj_T^z=\overline{\begin{pmatrix}\mathcal{T}_{1}(z)\\
\mathcal{T}_{2}(z)\\\cdot\\ \mathcal{T}_T(z)\end{pmatrix}}.
\end{equation}
\end{theorem}
Note \cite{MM6,MM8} that control $J_T j_T^z$ drives the system
(\ref{Jacobi_dyn}) to the special state at $t=T$:
\begin{equation*}
\left(W^TJ_Tj_T^z\right)_i=\left(W^T_TJ_Tj_T^z\right)_i=\overline{p_i(z)},\quad
i=1,\ldots,T.
\end{equation*}
Thus the reproducing kernel in $B^T_{A}$ is given by
\begin{multline*}
J_z^T(\lambda)=\sum_{k=1}^T \mathcal{T}_k(\lambda)(j^z_T)_{k}=\left(C_Tj^\lambda_T,j^z_T\right)_{\mathcal{F}^T}=\left(W^*_TW_Tj^\lambda_T,j^z_T\right)_{\mathcal{F}^T}\\
=\left(W_Tj_T^\lambda,W_Tj_T^z\right)_{\mathcal{H}^T}=\sum_{n=1}^T\overline{p_n(z)}{p_n(\lambda)}.
\end{multline*}

\begin{remark}
Due to the finite speed of wave propagation in (\ref{Jacobi_dyn}),
(\ref{Jacobi_dyn_int}) (cf. Remark \ref{Rem1}) we can use the
system (\ref{Jacobi_dyn_int}) with $N=T$ and use operator $C^T_T$
and the measure $d\rho_T(\lambda)$ in formulas for the scalar
product (\ref{DBr_Scal}).
\end{remark}

\subsection{Remark on canonical systems and Jacobi matrices.}

Let $2\times2$ matrix function $0\leqslant H=H^*\in
L_{1,loc}([0,L);\mathbb{R}^{2\times2})$ be a Hamiltonian,
define the matrix $J=\begin{pmatrix} &0&1 \\ -&1&0\end{pmatrix}$ and the vector $Y=\begin{pmatrix} Y_1 \\
Y_2\end{pmatrix}$ be solution to the following Cauchy problem:
\begin{equation}
\begin{cases}
\label{CanSyst}-J\frac{dY}{dx}=\lambda HY,\\
Y(0)=C,
\end{cases}
\end{equation}
for $C\in\mathbb{R}^2$, $C\not=0$. Without loss of generality we
assume that $\operatorname{tr}H(x)=1$.  Then the function
$E_x(\lambda)=Y_1(x,\lambda)+iY_2(x,\lambda)$ is a Hermite-Biehler
function ($E_L(\lambda)$ makes sense if $L<\infty$), it is called
de Branges function of the system (\ref{CanSyst}) since one can
construct de Branges space based on this function. On the other
hand, $E_L$ serves as an inverse spectral data for the canonical
system (\ref{CanSyst}). The main result of the theory
\cite{RR,DBr} says that every Hermite-Biehler function satisfying
some additional conditions comes from some canonical system.

\begin{remark}
The Jacobi matrices are particular examples of canonical system
(\ref{CanSyst}). The Hamiltonian $H(x)$ defined on the interval
$(0,L)$ related to the matrix $A$ is piecewise constant and has a
special structure \cite{RR}. Moreover, the Jacobi matrix $A$ is in
the limit circle case if and only if $L<+\infty$.
\end{remark}

This remark leads to the following question: is it possible to
introduce the de Branges spaces by the method from Section 5 for
semi-infinite Jacobi matrices without passing to canonical
systems? If the answer is ``yes'' then this method should ``feel'' the
difference between limit point and limit circle cases for Jacobi operator $A$, and
the norm in the space can be expressed in the dynamic terms as in
finite-dimensional situation.

\section{Properties of connecting operator and Hilbert spaces of functions associated with semi-infinite Jacobi matrices.}

\subsection{Operator $C$}

We consider the  matrix $C$ formally defined by the product of two
matrices: $C=\left(W\right)^* W$ (cf. (\ref{WT_1}),
(\ref{CT_def})), where
\begin{equation*}
W=\begin{pmatrix}
a_0 & w_{1,1} & w_{1,2} & \ldots \\
0 & a_0a_1 & w_{2,2} &  \ldots \\
\cdot & \cdot & \cdot & \cdot  \\
0 & \ldots & \prod_{j=1}^{k-1}a_j& \ldots \\
\cdot & \cdot & \cdot  & \cdot  \\
0 & 0 & 0  & \ldots
\end{pmatrix}
\end{equation*}
and
\begin{equation*}
C=
\begin{pmatrix}
r_0 & r_1& r_2 & \ldots & \ldots \\
r_1& r_0+r_2 & r_1+r_3 & \ldots &\ldots\\
r_2 & r_1+r_3 & \cdot & \cdot & \cdot \\
\ldots & \ldots & \ldots & \ldots & \ldots \\
\ldots & \ldots & \ldots & \ldots & \ldots
\end{pmatrix},\quad \{C\}_{ij}=\sum\limits_{k=0}^{\max{\{i,j\}}-1}r_{|i-j|+2k}.
\end{equation*}

In our approach to de Branges spaces the finite matrices $C_T$
plays an important role since they are used in the scalar product
(\ref{JM_scalprod}). We suggest that matrix $C $ should  play the
same role in semi-infinite case.


The matrix $C_T$ is connected with the classical Hankel matrix
$S_T$ by the following rule:
\begin{equation*}
C_T=\Lambda_T S_T\left(\Lambda_T\right)^*
\end{equation*}
\begin{proposition}
The entries of the matrix $\Lambda_T\in \mathbb{R}^{T\times T}$
are given by
\begin{equation*}
\Lambda_T=\left\{\alpha_{ij}\right\}=\begin{cases} 0,\quad \text{if $i<j$},\\
0,\quad \text{if $i+j$ is odd,}\\
D_{\frac{i+j}{2}}^j(-1)^{\frac{i+j}{2}+j},\quad \text{otherwise}
\end{cases}
\end{equation*}
where $D_n^k$ are binomial coefficients. The entries of the
response vector are related to moments by the rule:
\begin{equation}
\label{Resp_and_moments}
\begin{pmatrix}
T_1(\lambda)\\
\ldots \\
T_{T}(\lambda)
\end{pmatrix}=\Lambda_T\begin{pmatrix}
1\\
\lambda \\
\ldots \\
\lambda^{T-1}
\end{pmatrix},\quad
\begin{pmatrix}
r_0\\
r_1\\
\ldots \\
r_{T-1}
\end{pmatrix}=\Lambda_T\begin{pmatrix}
s_0\\
s_1\\
\ldots \\
s_{T-1}
\end{pmatrix}.
\end{equation}
\end{proposition}
We rewrite the Krein equation (\ref{Krein_eq}) in terms of $S_T$:
taking into account (\ref{Resp_and_moments}) we get that
\begin{equation*}
S_T\Lambda_T^*j_T^\lambda=\overline{\begin{pmatrix}
1\\
\lambda \\
\ldots \\
\lambda^{T-1}
\end{pmatrix}},
\end{equation*}

And introducing the vector $f_T^\lambda:=\Lambda_T^*j_T^\lambda$
we come to the equivalent form of (\ref{Krein_eq}):
\begin{equation}
\label{Krein_eq2}
S_T f_T^\lambda=\overline{\begin{pmatrix}
1\\
\lambda \\
\ldots \\
\lambda^{T-1}
\end{pmatrix}}.
\end{equation}
Therefore the reproducing kernel in $B^T_{A}$ is  represented by
\begin{equation*}
J^T_z(\lambda)\!=\!\left(C_T
j^\lambda_T,j^z_T\right)_{\mathcal{F}^T}=
\left(\overline{\begin{pmatrix}1
\\ \ldots\\ \lambda^{T-1}\end{pmatrix}},\Lambda_T^*j^z_T\right)_{\!\!\!\mathcal{F}^T}\!\!=\!\sum_{n=0}^{T-1}f_{T,k}^z\lambda^k.
\end{equation*}

Krein equations in the form (\ref{Krein_eq}) and in the form
(\ref{Krein_eq2}) demonstrate the importance of the knowledge of
the invertability properties of operators $S_T$ and $C_T$ when $T$
goes to infinity. The Theorem \ref{BergTH} answers this question
for $S$, below we answer the same question for $C$.

We introduce the notation:
\begin{equation*}
\beta_T:=\min\bigl\{\gamma_k\,|\,\gamma_k\text{ is an eigenvalue
of $C_T$},\, k=1,\ldots,T \bigr\}.
\end{equation*}
Then we can formulate the ``dynamic'' analog of Theorem
\ref{BergTH}.
\begin{theorem}
If the moment problem associated with sequence $\{s_k\}$ is
indeterminate $($the matrix $A$ is in the limit circle case$)$
then
\begin{equation}
\label{Beta0} \lim_{T\to\infty}\beta_T\geqslant
\left(\int\limits_{-1}^{1}l^{-1}(x)\frac{dx}{\sqrt{1-x^2}}\right)^{-1},\quad
l^{-1}(z)=\sum_{k=0}^\infty|p_k(z)|^2,
\end{equation}
\end{theorem}
\begin{proof}
We use the variational principal which says that
\begin{equation*}
\beta_T=\min\bigl\{\left(C_Tf,f\right)_{\mathcal{F}^T}\,|\,
\|f\|^2=1\bigr\}.
\end{equation*}
Passing to reciprocal gives
\begin{equation}
\label{Beta1}  \frac{1}{\beta_T}=\max\bigl\{\|f\|^2\,|\,
\left(C_Tf,f\right)_{\mathcal{F}^T}=1\bigr\}.
\end{equation}
Take functions $F,G\in B^T_{A}$, then
\begin{equation*}
F(\lambda)=\sum_{k=1}^T f_k\mathcal{T}_{k}(\lambda),\quad
G(\lambda)=\sum_{k=1}^T g_k\mathcal{T}_{k}(\lambda)
\end{equation*}
for some $f,g\in\mathcal{F}^T$. Since  $T_k(\lambda)$ are
Chebyshev polynomials of the second kind, they are orthogonal with
respect to the measure
\begin{equation*}
d\nu(\lambda)=\frac{\chi_{(-1,1)}(\lambda)}{\sqrt{1-\lambda^2}}d\lambda,
\end{equation*}
then
\begin{eqnarray*}
\|f\|^2=\sum_{k=1}^T|f_k|^2=\int\limits_{-1}^{1}|F(\lambda)|^2\frac{d\lambda}{\sqrt{1-\lambda^2}},\\
(C_Tf,g)_{\mathcal{F}^T}=\int\limits_{\mathbb{R}}\overline{F(\lambda)}{G(\lambda)}\,d\rho(\lambda),\quad
d\rho{\lambda}\in \mathcal{M}_H.
\end{eqnarray*}
Let us take another representation of $F(\lambda)$:
\begin{equation*}
F(\lambda)=\sum_{k=1}^T c_kp_k(\lambda)
\end{equation*}
Then for such $F$ we have that
\begin{equation*}
\int\limits_{\mathbb{R}}|F(\lambda)|^2\,d\rho(\lambda)=(C_Tf,f)_{\mathcal{F}^T}=\sum_{k=1}^T
|c_k|^2.
\end{equation*}
Thus we can rewrite (\ref{Beta1}) as:
\begin{equation}
\label{Beta2}
\frac{1}{\beta_T}=\max\left\{\int\limits_{-1}^{1}\overline{\sum_{i=1}^T
c_i p_i(\lambda)}{\sum_{j=1}^T c_j
p_j(\lambda)}\frac{d\lambda}{\sqrt{1-\lambda^2}}\,\Big|\,
\sum_{k=1}^T |c_k|^2=1\right\}
\end{equation}
Introducing the notation
$k_{ij}:=\int\limits_{-1}^{1}\overline{p_i(\lambda)}{p_j(\lambda)}\frac{d\lambda}{\sqrt{1-\lambda^2}}$,
we rewrite (\ref{Beta2}) as
\begin{equation}
\label{Beta3} \frac{1}{\beta_T}=\max\left\{\sum_{i,j=1}^T
k_{ij}\overline{c_i}{c_j}\,|\, \sum_{k=1}^T |c_k|^2=1\right\}.
\end{equation}
Define the matrix $K:=\left\{k_{ij}\right\}_{i,j=1}^T$, being a
Gram matrix, $K$ is positive.  The latter implies that the right
hand side of (\ref{Beta3}) is monotonically increasing, and
consequently, $\beta_T$ is monotonically decreasing as
$T\to\infty,$ which guarantees the existence of the limit in
(\ref{Beta0}). Then we proceed with the estimate:
\begin{equation*}
 \frac{1}{\beta_T}\leqslant
 \operatorname{Tr}{K}=\int\limits_{-1}^1\frac{\sum\limits_{k=1}^T|p_k(\lambda)|^2\,d\lambda}{\sqrt{1-\lambda^2}}\leqslant \int\limits_{-1}^1\frac{l^{-1}(\lambda)\,d\lambda}{\sqrt{1-\lambda^2}}.
\end{equation*}
The latter inequality yields the statement of the theorem.
\end{proof}

\begin{remark}
We note that the characterization limit point/limit circle in
terms of $\beta_T$ as in the Theorem \ref{BergTH} does not hold.
The simple example of free Jacobi operator, i.e. when $a_k=1$,
$b_k=0$, $k=1,2,\ldots$ confirms this. In this case it is not hard
to see that $C_T=I_T$ is an identity operator and consequntly,
$\beta_T=1$ for any $N$.
\end{remark}


We introduce the notation:
\begin{equation*}
\gamma_T=\max\bigl\{\gamma_k\,|\,\gamma_k\text{ is an eigenvalue
of $C_T$},\, k=1,\ldots,T \bigr\}.
\end{equation*}
\begin{lemma}
If  there exist such a constant $M\in\mathbb{R}$ that
$\gamma_T\leqslant M$ for all $T=1,2,\ldots$, then operator $A$ is
in the limit point case $($the moment problem associated with the
sequence $\{s_k\}$ is determined$)$.
\end{lemma}
\begin{proof}
The condition of the statement implies that the spectrum of
$C_T^{-1}$ is contained in $[\frac{1}{M},+\infty)$. Then we can
estimate the quadratic form
\begin{equation}
\label{Beta01}
\left(\left(S_T\right)^{-1}\xi,\xi\right)=\left(\left(C_T\right)^{-1}\Lambda_T\xi,\Lambda_T\xi\right)\geqslant
\frac{1}{M}\left(\Lambda_T\xi,\Lambda_T\xi\right),\quad \xi\in
\mathcal{F}^T.
\end{equation}
Choosing $\xi=e_1=(1,0,\ldots,0)$ we observe that due to
(\ref{Resp_and_moments}) we have that
\begin{equation}
\label{Beta02}
\Lambda_T\begin{pmatrix}1\\0\\ \cdot \\
0\end{pmatrix}=\begin{pmatrix}T_1(0)\\T_2(0)\\ \cdot \\
T_T(0)\end{pmatrix}.
\end{equation}
Since it is known that
\begin{equation}
\label{Beta03} T_{2n-1}(0)=(-1)^{n-1},\quad T_{2n}(0)=0,
\end{equation}
using (\ref{Beta01}), (\ref{Beta02}), (\ref{Beta03}), we come to
the estimate
\begin{equation*}
\left(\left(S_T\right)^{-1}e_1,e_1\right)\geqslant\frac{1}{M}\frac{T}{2}.
\end{equation*}
The latter inequality means that the maximal eigenvalue of
$S_T^{-1}$ tends to infinity when $T$ goes to infinity, which in
turn implies that the minimal eigenvalue of $S_T$ tends to zero
when $T$ goes to infinity. That yields the limit point case due to
Theorem \ref{BergTH}.
\end{proof}

The  matrix $C=\left\{c_{n,m}\right\}_{n,m=0}^\infty$ give rise to
the (formally defined) operator $\mathcal{C}$ in $l_2$:
\begin{equation*}
\left(\mathcal{C}f\right)_n:=\sum_{m=0}^\infty c_{nm}f_m,\quad
f\in l_2.
\end{equation*}
Then without any a'priory assumptions on $c_{nm}$, only the
quadratic form
\begin{equation*}
C[f,f]:=\sum_{m,n\geqslant 0}c_{nm}\overline{f_m}{f_n}
\end{equation*}
is well-defined on the domain $D$ (\ref{domA}).

We always assume the positivity condition
\begin{equation*}
\sum_{m,n\geqslant 0}c_{nm}\overline{f_m}f_n\geqslant 0,\quad f\in
D,
\end{equation*}
which guarantees \cite{MM10} the existence Jacobi matrix $A$ and
the (not necessarily uniquely defined) measure $M\in
\mathcal{M}_H$ which solves the moment problem
(\ref{MikhaylovAS_Moment_eq}). The following result is valid

\begin{theorem}
The quadratic form $C[,]$ is closable in $l_2$ if one of the
following occurs:
\begin{itemize}
\item{a)} If the Jacobi operator $A$ is in the limit circle case
$($the moment problem associated with sequence $\{s_k\}$ is
indetermined$)$ \item{b)} If operator $A$ is bounded, i.e. for all
$k=1,2,\ldots$, $|a_k|,|b_k|\leqslant c$ for some $c\in
\mathbb{R}_+$ and its spectral measure is absolutely continuous
with respect to Lebesgue measure.
\end{itemize}
\end{theorem}
\begin{proof}
Consider the operator $B: l^2\mapsto L_2(\mathbb{R},M)$, $M\in
\mathcal{M}_H$, acting by the rule
\begin{equation*}
Bf:=\sum_kf_k \mathcal{T}_k(\lambda),
\end{equation*}
with the domain $D$.
In view of (\ref{JM_scalprod}) the relationship between $B$ and
$C[,]$ is of the form
\begin{equation*}
C [f,f]=\|Bf\|^2_{L_2(\mathbb{R},M)},\quad f\in D.
\end{equation*}
By definition, the form $C[,]$ is closable in $l_2$ if and only if
the operator $B$ is closable.

We need to show that for the sequence $D\ni f^{(n)}\to 0$ in $l_2$
such that and $Bf^{(n)}\to F$ in $L_2(\mathbb{R},M)$ we
necessarily have that $F=0$.

In the limit circle case we use the existence of the holomorphic
kernel (\ref{Kenrnel0}), which gives that the limit function is
analytic: $F\in \mathcal{E}$. Indeed the sequence $f^{n}\in D$, so
$\mathcal{E} \ni Bf^{(n)}\to F$ in $L_2(\mathbb{R},M)$ then (see
the end of the Section 3) $F\in \mathcal{E}$. Alternatively, to
get the same result, we can pass to the limit in (\ref{Kernel}).
Since $f^{(n)}\to 0$ in $l_2,$ we have that corresponding
\begin{equation*}
Bf^{(n)}=:F^{(n)}\to 0\quad \text{in}\quad
L_2\left(-1,1;\frac{d\lambda}{\sqrt{1-\lambda^2}}\right),
\end{equation*}
The latter immediately yields $F=0$ in $L_2(-1,1;d\rho)$, so $F=0$
on $(-1,1)$ and thus being analytic in $\mathbb{C}$, $F=0$ in
$\mathbb{C}$. This implies the closability of the operator $B$ and
of the quadratic form $C[,]$.

In the case of bounded Jacobi matrix $A$ we can assume that
$\|A\|<1$ (otherwise considering the operator $\alpha A$ with
appropriate $\alpha$). In this case we have that the solution of
the moment problem is unique and denoted by $M$. The support of
theis measure is bounded: $\operatorname{supp}M\subset (-1,1)$.
Repeating the above arguments we get that for the limit function
$F$, which is in this case is not necessarily analytic, we also
have that $F=0$ in $L_2(-1,1;d\rho)$ and thus $F=0$ almost
everywhere on $(-1,1)$. But then we immediately obtain that $F=0$
in $L_2(\mathbb{R},M)$ since the support of $M\subset (-1,1)$ and
$M$ is absolutely continuous with respect to $dx$. And thus $C[,]$
and $B$ are closable.
\end{proof}

Now we assume the sequence $s_0,s_1,\ldots$ is indeterminate (the
matrix $A$ is in the limit circle case).  By analogy with
(\ref{De_Br_sp}) we introduce the linear manifold
\begin{equation*}
B_A^\infty:=\left\{\sum_{k=1}^\infty f_k \mathcal{T}_k(\lambda)\
\bigl{|}\ C[f,f]<\infty \right\}
\end{equation*}
The scalar product in $B_A^\infty$ is given by the rule
\begin{eqnarray*}
[F,G]_{B_A^\infty}:=C[f,g]_{}=\int\limits_{-\infty}^\infty
\overline{F(\lambda)}{G(\lambda)}\,d\,M, \quad M\in
\mathcal{M}_H,\\
F,G\in B_A^\infty,\quad F(\lambda)=\sum_{k=1}^\infty f_k
\mathcal{T}_k(\lambda),\quad G(\lambda)=\sum_{k=1}^\infty g_k
\mathcal{T}_k(\lambda).
\end{eqnarray*}
The reproducing kernel is given by (\ref{Kenrnel0}):
\begin{equation*}
J^\infty_z(\lambda)=\sum_{n=1}^\infty \overline{p_n(z)}{p_n(\lambda)}.
\end{equation*}
The conditions of Theorem \ref{TeorDB} verifying $B_A^\infty$ is a
de Branges space are trivially checked.


\begin{thebibliography}{99}

\bibitem{AH} {N.I. Akhiezer}. The classical moment problem.
\textit{Oliver and Boyd, Edinburgh.} 1965.




\bibitem{AT} {F. V. Atkinson}. Discrete and Continuous Boundary Problems.
\textit{New York/London, Academic Press.} 1964.





\bibitem{B07}
{M. I.Belishev}. Recent progress in the boundary control method.
\textit{Inverse Problems}, 23, no 5, R1--R67, 2007.


\bibitem{B17} {M.I.Belishev.} Boundary control and tomography of Riemannian
manifolds (the BC-method). \textit{Uspekhi Matem. Nauk,} 72, no.
4, 3--66 (in Russian), 2017.






\bibitem{Berg1}
{Berg C., Chen Y., Ismail M.E.H.} {Small eigenvalues of large
Hankel matrices: The indeterminate case.} \textit{Mathematica
Scandinavica,} 91, 67--81, 2002.

\bibitem{Berg2}
{Berg C., Szwarc R.} {Inverse of infinite hankel moment matrices}
\textit{Symmetry, Integrability and Geometry: Methods and
Applications (SIGMA),} 14, 2018,

\bibitem{Berg3}
{Berg C., Szwarc R.} {Closable Hankel Operators and Moment
Problems.}  \textit{Integral Equations and Operator Theory} 92,
no. 1,5, 2020.



\bibitem{DBr}
{Louis de Branges}. Hilbert space of entire functions.
\textit{Prentice-Hall, NJ.} 1968.



\bibitem{DMcK}
{H Dym, H. P. McKean}. Gaussian processes, function theory, and
the inverse spectral problem. \textit{Academic Press, New York
etc.} 1976.









\bibitem{MM6} {A.S. Mikhaylov, V.S. Mikhaylov.} {Inverse dynamic problems for
canonical systems  and de Branges spaces.} \textit{Nanosystems:
Physics, Chemistry, Mathematics,} 9, no. 2, 215-224, 2018.

\bibitem{MM7} {A.S. Mikhaylov, V.S. Mikhaylov.} {Boundary Control
method and de Branges spaces. Schr\"odinger operator, Dirac
system, discrete Schr\"odinger operator.} \textit{Journal of
Mathematical Analysis and Applications,} 460, no. 2, 927-953,
2018.

\bibitem{MM8} {A. S. Mikhaylov, V. S. Mikhaylov.} {Dynamic inverse problem for
Jacobi matrices.} \textit{Inverse Problems and Imaging,} 13, no.
3, 431-447, 2019.

\bibitem{MM9} {A. S. Mikhaylov, V. S. Mikhaylov.} {On application of the
Boundary control method to classical moment problems.}
\textit{Journal of Physics: Conference Series,} 2092(1):012002,
2021.


\bibitem{MM10} {A. S. Mikhaylov, V. S. Mikhaylov.} {Inverse problem for
dynamical system associated with Jacobi matrices and classical
moment problems.} \textit{Journal of Mathematical Analysis and
Applications,} 487, no. 1, 2020.

\bibitem{MM11}
{A.S. Mikhaylov, V.S. Mikhaylov,} On a special Hilbert space of
functions associated with the multidimensional wave equation in a
bounded domain.  \textit{IEEE Conference Proceedings: Days on
Diffraction'2022.} 101 -- 105, 2022.


\bibitem{R2} {C. Remling }. Schr\"odinger operators and de Branges spaces. \textit{J.
Funct. Anal.} 196, no. 2, 323--394, 2002.

\bibitem{RR} {R. V. Romanov.} Canonical systems and de Branges spaces
\textit{http://arxiv.org/abs/1408.6022}.

\bibitem{S}
{B. Simon}, {The classical moment problem as a self-adjoint finite
difference operator.},  \textit{Advances in Math.,} 137, 82-203,
1998.

\bibitem{KSch} {K. Schm\"udgen} {The moment problem.}
\textit{Springer International Publishing, Cham,} 2017.

\bibitem{Y}
{D.R. Yafaev} {Unbounded Hankel Operators and Moment Problems},
\textit{Integral Equations and Operator Theory,}  85, no. 2,
289--301, 2016.

\bibitem{Y2}
{D.R. Yafaev} Correction to: Unbounded Hankel Operators and Moment
Problems (Integral Equations and Operator Theory, (2016), 85, 2,
(289-300),) \textit{Integral Equations and Operator Theory,} 91,
no. 5, 2019.


\end{thebibliography}
\end{document}